\documentclass{svjour3}                     
\smartqed  
\usepackage{graphicx}
 \usepackage{color}
\usepackage{amsmath}
\usepackage{amssymb}
\usepackage[latin1]{inputenc}
\usepackage{theorem}
\usepackage{natbib}
\usepackage{epsfig}

\newcommand{\ds}{\displaystyle}
\def\text#1{\quad\mbox{#1}\quad} 
\def\mtext#1{\,\mbox{#1}\,} 
\def\defegal{:=}

\def\NN{{\mathbb N}}
\def\RR{{\mathbb R}}
\def\UU{{\mathbb U}}
\def\XX{{\mathbb X}}

\def\obj{{\mathbb D}}
\newcommand{\VV}{\mathbb{V}}

\def\SCS{\VV_0}

\def\SSB{\mbox{$S\!S\!B$}}
\def\spr{\mbox{spr}}
\def\stock{{\rm stock}}

\newcommand{\mmin}{{\sf min}}
\def\yield{{\sf yield}}
\def\eco{{\sf protect}}

\def\Yield{Y}

\newcommand{\llim}{{\sf lim}}
\newcommand{\Blim}{B_{\llim}}
\newcommand{\Flim}{F_{\llim}}

\def\biomass{B}
\def\dynamics{g}
\def\catch{C}
\def\multiplier{\lambda}
\def\abundance{N}
\def\mature{\gamma}
\def\weight{w}
\def\exploitation{F}
\def\mortality{M}
\def\ind#1{_{#1}}
\def\spawning_stock_biomass{{\sf SSB}}
\def\fishing_mortality{{\sf F}}
\def\stock_recruitment{{\sf S-R}}

\newcommand{\upper}{_{\sharp}}
\newcommand{\llower}{_{\flat}}

\def\control{u} 
\def\state{x} 

\def\constraint{{\cal L}}
\def\threshold{l}

\renewcommand{\omega}{w}

\def\dynamics{G}

\def\1{{\mathbf 1}}
\def\upper{^{\sharp}}
\def\llower{^{\flat}}

\def\defegal{:=}

\journalname{J. Math. Biol.}
\begin{document}

\title{Viable harvest of monotone bioeconomic models}

\author{M. De Lara        \and P. Gajardo \and H.  Ramírez C.}

\institute{M. De Lara  \at Universit\'e Paris-Est, Cermics,
6-8 avenue Blaise Pascal, 77455 Marne la Vall\'ee Cedex 2, France.
             \\
              \email{delara@cermics.enpc.fr}                    \and
           P. Gajardo  \at
             Departamento de Matemática, 
Universidad Técnica Federico Santa María, 
Avda. España 1680 Casilla 110-V, Valparaíso, Chile. \\
              \email{pedro.gajardo@usm.cl}                    \and
           H. Ramírez C.  \at
           Universidad de Chile, Centro de Modelamiento Matem\'atico
  (CNRS UMI 2807) and Departamento de Ingenier\'{\i}a Matem\'atica, 
	Casilla 170/3, Correo 3, Santiago, Chile.\\
              \email{hramirez@dim.uchile.cl}         }

\date{Received: date / Accepted: date}

\maketitle

\begin{abstract}
Some monospecies age class models, as well as specific multi-species
models (with so-called technical interactions), exhibit useful 
monotonicity properties.
This paper deals with discrete time monotone bioeconomic dynamics in
the presence of state and control constraints.
In practice, these latter ``acceptable configurations''
represent production and preservation
requirements to be satisfied for all time, 
and they also possess monotonicity properties. 
A state $\state$ is said to belong to the viability kernel 
if there exists a trajectory, of states and controls, starting from
$\state$ and satisfying the constraints.
Under monotonicity assumptions, we present upper and lower estimates of 
the viability kernel.  
This helps delineating domains where a viable management is possible.
Numerical examples, in  the context of fisheries management, 
for the Chilean sea bass (\emph{Dissostichus eleginoides}) 
and Alfonsino (\emph{Beryx splendens})  are   given.

\end{abstract}

\section{Introduction}

This paper deals with the  control of discrete--time dynamical systems 
of the form\\ $\state(t+1)=\dynamics\big(\state(t),\control(t)\big)$, $t \in \NN$, 
with state $\state(t) \in \XX$ and control $\control(t) \in \UU$,
in the presence of state and control constraints $\big(\state(t),\control(t)\big)
\in \obj$.  The subset $\obj \subset \XX \times \UU$
describes ``acceptable configurations of the system''.
Such problems of dynamic control under constraints refers to viability 
\cite{Aubin:1991} or invariance \cite{Clarkeetal:1995} framework.
From the mathematical viewpoint, most of viability and weak invariance
results are addressed in the continuous time case. However, some
mathematical works deal with the discrete-time case. This  includes the
study of numerical schemes for the approximation of the viability
problems of the continuous dynamics as
in \cite{Aubin:1991,Stpierre:1994,QSP:1995}.  
In the control theory literature, problems of constrained control have
also been addressed in the discrete time case (see the survey paper
\cite{Blanchini:1999});
reachability of target sets or tubes for nonlinear discrete time
dynamics is examined in \cite{bertsekas71minimax}. 

We consider sustainable management issues which can be formulated within
such a framework  as in
\cite{benedoyen:2000,doyenbene2003,Bene-Doyen-Gabay:2001,Eisenack-Sheffran-Kropp:2006,martinetdoyenx,Mullonetal:2004,Rapaport-Terreaux-Doyen:2006,DeLara-Doyen-Guilbaud-Rochet_IJMS:2007,DeLara-Doyen:2008}.

The time index $t$ is an integer and the time period $[t,t+1[$ may be a
year, a month, etc. 
The dynamics is generally a population dynamics, with state vector $\state(t)$
being either the biomass of a single species,
or a couple of biomasses for a predator--prey system,
or a vector of abundances at ages for one or for several species,
or abundances at different spatial patches, etc.
The control $\control(t)$ may represent catches, fishing mortality or fishing effort.
The ``acceptable set'' $\obj$ such that 
$\big(\state(t),\control(t)\big) \in \obj$
may include biological, ecological and economic
objectives as in \cite{Bene-Doyen-Gabay:2001}.
For instance, if the state $x$ is a vector of abundances at ages
and the control $u$ is a fishing effort, 
$\obj = \{ (x,u) \mid B(x) \geq b\llower \; , \; 
E(x,u) \geq e\llower \}$ 
represents acceptable configurations where 
conservation is ensured by an biological indicator $B(x) \geq b\llower$
(spawning stock biomass above a reference point, for instance)
and economics is taken into account via minimal catches 
$E(x,u) \geq e\llower$ (catches $E(x,u)$ above a threshold).

The viability kernel $\VV(\dynamics,\obj )$ associated with the dynamics
$\dynamics$ and the acceptable set\footnote{%
In \cite{Aubin:1991}, the viability kernel $\VV_K(\dynamics)$ is defined
with respect to the dynamics $\dynamics$ and to a subset $K \subset \XX$
of the state space $\XX$, and the constraints on the controls are 
contained in the definition in $\dynamics$.
We prefer to put together the set of state constraints with the set of
admissible controls, although these sets play very different roles. 
Indeed, in practice, constraints are expressed \emph{via} indicators 
which are functions of both variables (state and control), 
especially for production constraints which depend on the catches.
Thus, the set $\obj$ makes the conflicting requirements, 
between preservation and production, more visible than with the 
Aubin's formalism.
}
 $\obj$ is known to play a basic role
for the analysis of such problems and the design of viable control
feedbacks. Unfortunately, its computation is not an easy task in
general. However, following an approach initiated
in~\cite{DeLara-Doyen-Guilbaud-Rochet_SCL:2006}, the viability kernel
may be estimated from below or form above under specific monotonicity
assumptions which are pertinent for a class of bioeconomic models. 

In Sect.~\ref{sec:Viability_issues_and_monotone_bioeconomic_models}, 
we recall the viability issues in discrete time, 
and we introduce monotone bioeconomic models.
Sect.~\ref{sec:Viability_kernel_estimates_for_monotone_bioeconomic_models} 
provides our main theoretical results on estimates of the viability kernel.
An application to fishery management is provided in
Sect.~\ref{sec:applications} with numerical estimates for the 
Chilean sea bass
(\emph{Dissostichus eleginoides}), harvested in the south of Chile, and
Alfonsino (\emph{Beryx splendens}), harvested in the Juan Fernández
archipelago.

\section{Viability issues and monotone bioeconomic models}
\label{sec:Viability_issues_and_monotone_bioeconomic_models}

In this introductory section, we recall the viability issues in discrete time, 
then introduce monotone bioeconomic models.

\subsection{Viability in discrete time}

Let us  consider a nonlinear control system  described in discrete
time by the difference equation
\begin{equation*}
\left\{ \begin{array}{l}
\state(t+1)=\dynamics(\state(t),\control(t)), \quad
t=t_{0},t_{0}+1,\ldots \\
\state(t_{0})\quad\mtext{given,}
\end{array} \right .
\label{eq:generaldyn}
\end{equation*}
where  the \emph{state variable} $\state(t)$ belongs to  the
finite dimensional state space $\XX \subset \RR^{n_{\XX}}$, the
\emph{control variable} $\control(t)$ is an element of the
\emph{control set} $\UU \subset \RR^{n_{\UU}}$ while the \emph{dynamics}
$\dynamics$ maps $\XX \times \UU$ into  $\XX$. In our context,
$\state(t)$ will typically represent the vector of abundances per
age class of a population, while $\control(t)$ will be a harvest
(induced mortality, harvesting effort, etc.).

A decision maker  describes \emph{acceptable configurations of the
system} through a set $\obj  \subset \XX \times \UU$ termed the
\emph{acceptable set}
\begin{equation*}
    (\state(t),\control(t)) \in \obj \, , \quad
    \forall t=t_{0},t_{0}+1,\ldots
\end{equation*}
where $\obj $ includes both system states and controls
constraints. Typical instances of such an acceptable set are given
by inequalities requirements
\begin{equation}
 \obj = \{(\state,\control) \in \XX
\times \UU \mid \forall i=1,\ldots,p \, , \quad
\constraint_i(\state,\control) \geq \threshold_i \},
\label{eq:constraints} 
\end{equation}
 where the
functions $\constraint_1$, \ldots, $\constraint_p$ may be
interpreted as {\em indicators}, and the real numbers $\threshold_1$,
\ldots, $\threshold_p$ as the corresponding {\em thresholds}. For
management issues, the set $\obj $ will be the mathematical
expression of preservation and/or production objectives.
\bigskip

Viability is defined 
as the ability to choose, at each time step
$t=t_{0},t_{0}+1,\ldots$, a control $\control(t) \in \UU$ such
that the system configuration remains acceptable. More precisely,
the system is viable  if the following feasible set is not empty:
\begin{equation*}
\VV(\dynamics,\obj ) \defegal \left\{\state \in \XX \left|
\begin{array}{l}\exists\; (\control(t_{0}),\control(t_{0}+1), \ldots )
\mtext{ and } (\state(t_{0}),\state(t_{0}+1), \ldots )\\[2mm]
\mtext{ satisfying } \state(t_{0})=\state \, , \quad
\state(t+1)=\dynamics(\state(t),\control(t)) \\[2mm]
\mbox{ and } (\state(t),\control(t)) \in \obj \, , \quad
    \forall t=t_{0},t_{0}+1,\ldots
\end{array} \right.
\right\} \; .
\label{eq:viability_kernel}
\end{equation*}
For a decision maker, knowing the viability kernel has practical
interest since it describes the set of states from which controls can be
found that maintain the system in an acceptable configuration
forever. However, computing this kernel is not an easy task in general. 

We shall focus on   estimates of viability
kernels when the dynamics $\dynamics$ and the acceptable sets have
specificic monotony properties.
For this purpose, we shall introduce a generic form for dynamics and
acceptable sets corresponding to what we shall call 
\emph{monotone population models harvesting issues}.

\subsection{Monotone bioeconomic models}

In what follows, the state space $\XX$ and the control space $\UU$
are subsets $\XX \subset \RR^{n_{\XX}}$ and $\UU \subset
\RR^{n_{\UU}}$ supplied with the componentwise order: $x' \geq x$
if and only if each component of $x'=(x'_1,\ldots,x'_{n_{\XX}})$ is
greater than or equal to the corresponding component of
$x=(x_1,\ldots,x_{n_{\XX}})$:
\( x' \geq x \iff x'_i \geq x_i, \;
i=1,\ldots,n_{\XX} \). 
 A mapping $f: \RR^a \to \RR^b$ is said to be
increasing if $x \geq x' \Rightarrow f(x)  \geq f(x')$. Similar
definition holds for decreasing.

\subsubsection*{Dynamics}

Monospecies dynamic population models generally have the
following qualitative properties:
i) the higher the state abundance vector, the higher at next period;
ii) the higher the harvest, the lower the state abundance vector at next
period. 
Some specific multi-species models, without ecological but with so
called technical interactions, share such properties.
This motivates the following definitions.

 We say that the dynamics $\dynamics: \; \XX\times
\UU \rightarrow \XX$ is \emph{increasing with respect to the
state} if it satisfies $ \forall~ (\state,\state',\control) \in
\XX \times \XX \times \UU $,
\( \state' \geq \state \Rightarrow
\dynamics(\state',\control) \geq \dynamics(\state,\control),\) 
and is \emph{decreasing with respect to the control} if 
\( \forall
(\state,\control,\control') \in \XX \times \UU \times \UU \) , 
\( \control' \geq \control \Rightarrow \dynamics(\state,\control')
\leq \dynamics(\state,\control) \).

We shall coin $\dynamics: \XX \times \UU \to \XX $ a
\emph{monotone bioeconomic dynamics} if $\dynamics$ is increasing
with respect to the state and decreasing with respect to the
control. 

\subsubsection*{Bounded control set}

Assuming they exist, we denote by $\control\llower,\,
\control\upper \in \UU$ the lower and upper bounds of the set
$\UU$, \emph{i. e.}  $\control\llower \le \control \leq
\control\upper$ for all $\control \in \UU$.

\subsubsection*{Upper and lower dynamics without control}

Let the dynamics $\dynamics$ be a monotone bioeconomic dynamics.
Define the \emph{upper\footnote{%
Because $\dynamics \leq \dynamics\llower$. } dynamics without
control} $\dynamics\llower:\XX\to \XX$ by $
\dynamics\llower(\state)=\dynamics (\state,\control\llower)$. Its
$t$ iterate ($t=t_0,t_0+1,\ldots$) will be denoted by $
(\dynamics\llower)^{(t)}$. In the same way, the \emph{lower
dynamics without control} is defined by $
\dynamics\upper(\state)=\dynamics (\state,\control\upper)$. With
these notations, we have that
\begin{equation*}
\dynamics\upper(\state) \leq \dynamics(\state,\control)  \leq
\dynamics\llower(\state) \, , \quad \forall ~(\state,\control) \in
\XX \times \UU \; . \label{eq:dynamics_bounds}
\end{equation*}

\subsubsection*{Acceptable set}

We say that a set $S \subset \XX$ is an \emph{upper set}
 (or is an \emph{increasing set})
if it satisfies the following property: $\forall \state \in S$, 
$\forall \state' \in \XX $ , $ \state' \geq \state
\Rightarrow \state' \in S$. In the same way, a set $K \subset
\XX\times \UU$ is said to be an \emph{upper set} if $\forall
(\state,\control) \in K $, $ \forall \state' \in \XX $ , $ \state'
\geq \state \Rightarrow (\state',\control) \in K$.

An acceptable set $\obj$ is said to be \emph{a production acceptable
set} if $\obj$ is
 \emph{increasing with respect both to the state and to the control}, that is
$\forall ~(\state,\state',\control,\control') \in
 \XX \times \XX \times\UU \times \UU $,
\( \state' \geq \state,~ \control' \geq \control,~~~
(\state,\control) \in \obj \Rightarrow (\state',\control') \in
\obj.\)
Particular instances are given by acceptable sets of the
form~\eqref{eq:constraints} 
 where the
indicators $\constraint_1$, \ldots, $\constraint_p$ are
increasing with respect to both variables (state and control). 
For instance, requiring a minimum yield may be captured by the
acceptable set
\( \obj_{\yield} =  \{(\state,\control) \mid
\Yield(\state,\control ) \geq y\llower \} \)
 where $y\llower \in \RR$ is a minimum yield threshold and
where the yield function $\Yield: \XX \times \UU \to  \RR$
is increasing with respect to both variables (state and control).

An acceptable set $\obj$ is said to be \emph{a preservation
acceptable set} if $\obj$ is
 \emph{increasing with respect to the state and
decreasing with respect to the control}, that is $\forall~
(\state,\state',\control,\control') \in
 \XX \times \XX \times\UU \times \UU $, 
 \(
 \state' \geq \state,~ \control' \leq \control,~~
(\state,\control) \in \obj \Rightarrow (\state',\control') \in
\obj. \)
Particular instances are given by acceptable sets of the
form~\eqref{eq:constraints} where the
indicators $\constraint_1$, \ldots, $\constraint_p$ are
increasing with respect to the state
but decreasing with respect to the control.
For instance, the ICES\footnote{%
International Council for the Exploration of the Sea.}
precautionary approach may be stated in the
viability framework 
with the following preservation acceptable set
\( \obj_{\eco} = \{(\state,\control)\in \XX \times \UU \mid
\SSB(\state) \geq B_{lim},\, F(\control) \leq \Flim \} \)
as in \cite{DeLara-Doyen-Guilbaud-Rochet_IJMS:2007}.
Here, $\SSB(\state)$ is the spawning stock biomass,
increasing with respect to the state,
while the fishing mortality $F(\control)$ is 
increasing\footnote{%
Hence $-F(\control)$ is decreasing with respect to the control.
To be consistent with the notation in~\eqref{eq:constraints}, it suffice
to rewrite \( \obj_{\eco} = \{(\state,\control)\in \XX \times \UU \mid
\SSB(\state) \geq B_{lim},\, -F(\control) \geq -\Flim \} \).
}
with respect to the control.

Notice that both production and preservation acceptable sets are
upper sets.
\bigskip

For any acceptable set $\obj$, introduce the 
\emph{state constraints set}
\begin{equation}
\SCS \defegal {\rm Proj}_{\XX}(\obj ) = \{\state \in \XX \mid
\exists \control \in \UU \, , \, (\state,\control) \in \obj \} \;
, \label{eq:SCS}
\end{equation}
obtained by projecting the acceptable set $\obj$ onto the state
space $\XX$. Introduce also
\begin{equation}
\left\{ \begin{array}{rclrcl}
 \SCS\upper &\defegal& \{\state \in \XX \mid
(\state,\control\upper) \in \obj \} \subset \SCS \; , \quad
& \obj\upper &\defegal & \SCS\upper \times \{ \control\upper \} \\[4mm]
\SCS\llower &\defegal& \{\state \in \XX \mid
(\state,\control\llower) \in \obj \} \subset \SCS \; , \quad &
\obj\llower &\defegal&  \SCS\llower \times \{ \control\llower \}
\; .
\end{array} \right.
\label{eq:SCSs}
\end{equation}
Notice that if  $\obj$ is a production acceptable set, we have
$\SCS = \SCS\upper$, and if $\obj$ is a preservation acceptable
set, we have $\SCS = \SCS\llower$.


\section{Viability kernel estimates for monotone bioeconomic models}
\label{sec:Viability_kernel_estimates_for_monotone_bioeconomic_models}


In this part, we shall provide lower and upper estimates of the
viability kernel $\VV(\dynamics,\obj )$ thanks to the following
sets $\VV(\dynamics\llower,\obj\llower)$,
$\VV(\dynamics\llower,\obj\upper)$, $
\VV(\dynamics\upper,\obj\llower)$ and
$\VV(\dynamics\upper,\obj\upper)$. These latter sets are easier to
compute than the viability kernel $\VV(\dynamics,\obj )$ 
because the dynamics
$\dynamics\llower$ and $\dynamics\upper$ have no control. Indeed,
by~\eqref{eq:SCS}, one obtains that, for any acceptable set $\obj$,
\begin{equation}
\begin{array}{rcl} \VV(\dynamics\llower,\obj ) &=&
\bigcap \limits_{t=0}^{+\infty} \{ \state \in \XX  \mid
(\dynamics\llower)^{(t)}(\state) \in \SCS \} \\[5mm] 
\VV(\dynamics\upper,\obj ) &=&
\bigcap \limits_{t=0}^{+\infty} \{ \state \in \XX  \mid
(\dynamics\upper)^{(t)}(\state) \in \SCS \} \; .\end{array} 
\label{eq:viability_kernels_no_control}
\end{equation}

\begin{proposition}\label{pr:viab-preservation-production}
Suppose that $\dynamics$ is a monotone bioeconomic dynamics and
that the control set $\UU$ has lower and upper bounds
$\control\llower,\, \control\upper \in \UU$.
\begin{enumerate}
\item If $\obj$ is a production acceptable set, then
\begin{equation}
\VV(\dynamics\upper,\obj\llower) \subseteq
\VV(\dynamics\llower,\obj\llower) \cup
\VV(\dynamics\upper,\obj\upper)
\subseteq \VV(\dynamics,\obj )\subseteq
\VV(\dynamics\llower,\obj\upper)
\; . \label{eq:VV-prod}
\end{equation}

\item  If $\obj$ is a preservation acceptable set, then
\begin{equation}
\VV(\dynamics\upper,\obj\upper) \subseteq \VV(\dynamics,\obj )=
\VV(\dynamics\llower,\obj\llower) \; .
\label{eq:VV-protect}
\end{equation}
\end{enumerate}

\end{proposition}

\begin{proof}
First, let us notice that, whatever the acceptable set $\obj$ and
the dynamics $\dynamics$, we have the inclusion
\begin{equation}
\VV(\dynamics\llower,\obj\llower) \cup
\VV(\dynamics\upper,\obj\upper) \subseteq \VV(\dynamics,\obj ) \;
. \label{eq:inclusion_1}
\end{equation}
Indeed,
\(
\VV(\dynamics\llower,\obj\llower) = \bigcap \limits_{t=0}^{+\infty}
\{ \state \in \XX  \mid \big( (\dynamics\llower)^{(t)}(\state) ,
\control\llower \big) \in \obj\} \subseteq \VV(\dynamics,\obj) \;
,
\)
since $\state \in \VV(\dynamics\llower,\obj\llower)$ means that
the stationary control $\control(t)=\control\llower$ makes that
the trajectory
$(\state(t),\control(t))=((\dynamics\llower)^{(t)}(\state),\control\llower)$
belongs to $\obj$. The same may be done with the control
$\control\upper$.

Second, when $\dynamics $ is a monotone bioeconomic dynamics and
$\obj$ is either a production or a preservation acceptable set,
we have the inclusions
\begin{equation}
\VV(\dynamics\upper,\obj ) \subset \VV(\dynamics,\obj ) \subset
\VV(\dynamics\llower,\obj ) \; . \label{eq:inclusion_2}
\end{equation}
This is a straightforward application of Proposition~11 in
\cite{DeLara-Doyen-Guilbaud-Rochet_SCL:2006}, because both $
\dynamics\llower$ and $\dynamics $ are increasing with respect to
the state, and because $\obj$ is an upper set (production and
preservation acceptable sets are upper sets).
\bigskip

Now, we come to the proof.
\begin{enumerate}

\item On the one hand, we have that
$\VV(\dynamics\upper,\obj\llower) \subseteq
\VV(\dynamics\llower,\obj\llower) $ by~\eqref{eq:inclusion_2} with
$\obj$ replaced by $\obj\llower$. By~\eqref{eq:inclusion_1}, this
gives the two lower estimates of the viability kernel
$\VV(\dynamics,\obj )$ in~\eqref{eq:VV-prod}.

On the other hand, since $\obj$ is a production acceptable set, we
have $ \SCS = \SCS\upper $, and thus, by~\eqref{eq:SCSs} and
\eqref{eq:viability_kernels_no_control},
\begin{equation*}
\begin{split}
 \VV(\dynamics\llower,\obj ) =
\bigcap_{t=0}^{+\infty} \{ \state \in \XX  \mid
(\dynamics\llower)^{(t)}(\state) \in \SCS \} = \bigcap_{t=0}^{+\infty}
\{ \state \in \XX  \mid (\dynamics\llower)^{(t)}(\state)
\in \SCS\upper \} \\[2mm] 
= \VV(\dynamics\llower,\obj\upper) \; .
\end{split}
\end{equation*}
As we have seen by~\eqref{eq:inclusion_2} that $\VV(\dynamics,\obj
) \subseteq \VV(\dynamics\llower,\obj ) $, this gives
$\VV(\dynamics,\obj ) \subseteq \VV(\dynamics\llower,\obj\upper)
$, hence the upper estimate of the viability kernel
$\VV(\dynamics,\obj )$ in~\eqref{eq:VV-prod}.

\item The lower estimate of the viability kernel
$\VV(\dynamics,\obj )$ in~\eqref{eq:VV-protect} comes
from~\eqref{eq:inclusion_1}.

Now, let us prove the equality  $
\VV(\dynamics\llower,\obj\llower) = \VV(\dynamics,\obj )$
in~\eqref{eq:VV-protect}. On the one hand,
by~\eqref{eq:inclusion_1} we know that $
\VV(\dynamics\llower,\obj\llower) \subseteq \VV(\dynamics,\obj )$.
On the other hand, since $\obj$ is a preservation acceptable set,
we have $ \SCS = \SCS\llower $, and thus
\begin{equation*}
\begin{split}
 \VV(\dynamics\llower,\obj ) =
\bigcap_{t=0}^{+\infty} \{ \state \in \XX  \mid
(\dynamics\llower)^{(t)}(\state) \in \SCS \} = \bigcap_{t=0}^{+\infty}
 \{ \state \in \XX  \mid (\dynamics\llower)^{(t)}(\state)
\in \SCS\llower \}\\[2mm] 
= \VV(\dynamics\llower,\obj\llower) \; .
\end{split}
\end{equation*}
By~\eqref{eq:inclusion_2}, this gives
$\VV(\dynamics\llower,\obj\llower) \subseteq \VV(\dynamics,\obj )
\subseteq \VV(\dynamics\llower,\obj )
=\VV(\dynamics\llower,\obj\llower) $.

\end{enumerate}

\end{proof}

When the acceptable set is given by means of indicators functions and
thresholds as in~\eqref{eq:constraints}, and the upper dynamics
$\dynamics\llower$ has a 
steady state satisfying some requirements,
 we obtain the following practical conditions for 
nonemptyness of the viability kernel. 

\begin{corollary}
\label{corollary1}
Suppose that $\dynamics$ is a monotone bioeconomic dynamics, that
the control set $\UU$ has lower and upper bounds $\control\llower,\,
\control\upper \in \UU$ and 
that the acceptable set $\obj$ is given by~\eqref{eq:constraints}
with indicators $\constraint_i$'s being upper semi-continuous functions
in the first (state) variable.  
Assume also that the upper dynamics $\dynamics\llower$ has a
steady state $\overline{\state}(\control\llower)$ 
 and there exists $L < 1$ such that
 \begin{equation}\label{condition_steady}
 \|\dynamics\llower(\state)-\overline{\state}(\control\llower)\| \le L \|\state - \overline{\state}(\control\llower)\| \qquad \forall \state \in \SCS
 \end{equation}
 for some norm $\|\cdot \|$ in $\XX$.

\begin{enumerate}
\item \label{it:cor_production}
If $\obj$ is a production acceptable set, one has
\begin{equation}
\exists i=1,\ldots,p \; , \quad 
\constraint_i(\overline{\state}(\control\llower),\control\upper) 
< \threshold_i \Rightarrow  \VV(\dynamics,\obj ) = \emptyset \; . 
\label{eq:VV-prod-cor}
\end{equation}

\item \label{it:cor_preservation}
  If $\obj$ is a preservation acceptable set, one has
\begin{equation}
\VV(\dynamics,\obj )\neq \emptyset \Leftrightarrow  \constraint_i(\overline{\state}(\control\llower),\control\llower) \ge \threshold_i \qquad  \forall  i=1,\ldots,p\;.
\label{eq:VV-protect-cor}
\end{equation}
\end{enumerate}
\end{corollary}
\begin{proof}
We proceed to prove Statement~\ref{it:cor_production} by
contra-reciprocal argument. 
Let us suppose that $\VV(\dynamics,\obj ) \neq \emptyset$ and take $\state$ in this set which is included in $\SCS$. From  Proposition~\ref{pr:viab-preservation-production}, $\state$ belongs to $\VV(\dynamics\llower,\obj\upper)$ or, equivalently
\[
 (\dynamics\llower)^{(t)}(\state)
 \in \SCS\upper=\SCS \; , \quad \forall t \geq t_{0} \Leftrightarrow 
 \constraint_i\big( (\dynamics\llower)^{(t)}(\state),\control\upper \big) \ge
 \threshold_i ~~~ \forall  i=1,\ldots,p   ~~~ \forall t \geq t_{0} \; .
\]
 Since $(\dynamics\llower)^{(t)}(\state) \in \SCS$ for all $t \geq
 t_{0}$, condition \eqref{condition_steady} implies
$(\dynamics\llower)^{(t)}(\state) \to
\overline{\state}(\control\llower)$.  Then,  from  the upper semi-continuity property of  functions $\constraint_i(\cdot, \control\upper)$, we obtain  the desired inequalities
\[
 \constraint_i(\overline{\state}(\control\llower),\control\upper) 
\ge \threshold_i \; , \quad  \forall  i=1,\ldots,p\; .
\]
The proof of necessary condition ($\Rightarrow$) in Statement~\ref{it:cor_preservation}  is analogous. For the sufficient condition ($\Leftarrow$) directly one can prove that $\overline{\state}(\control\llower) \in \VV(\dynamics,\obj)$ taking the stationary control $\control(t)=\control\llower$.
\end{proof}

The previous corollary provides necessary conditions (in the case of
a production acceptable set) and necessary and sufficient conditions (in
the case of a preservation acceptable set) to assure the non-emptiness of the
viability kernel. The quantities
$\constraint_i(\overline{\state}(\control\llower),\control\upper)$ (for
production acceptable sets) and
$\constraint_i(\overline{\state}(\control\llower),\control\llower)$ (for
preservation acceptable sets) can be interpreted as 
\emph{maximal thresholds} for the acceptable configurations. That is,
no trajectory $\state(\cdot)=(\state(t_0), \state(t_0 +1),\hdots)$ 
can generate values $\constraint_i(\state(t),\control(t))$ above these
values for all periods of time $t$, whatever the initial state
$\state_0$ and the control trajectory 
$\control(\cdot)=(\control(t_0), \control(t_0 +1),\hdots)$ be.

\begin{remark}
An alternative (non equivalent) condition to \eqref{condition_steady} in
the above Corollary~\ref{corollary1}
 would be supposing that the steady state
$\overline{\state}(\control\llower)$ is globally asymptotic stable on
$\SCS$ for the dynamics $\dynamics\llower$.
However, this is a strong assumption.
A weaker one would restrict global asymptotic stability on the subset 
$\VV(\dynamics,\obj) \subset \SCS$ 
(see the proof of Corollary~\ref{corollary1}). 
Nevertheless, it is not elegant neither practical to make any assumption
on the viability kernel $\VV(\dynamics,\obj)$, which is an object of
study and which might be empty.
\end{remark}

\section{Application to fishery management}
\label{sec:applications}

Now, we apply and specify the previous results in the case of an
age-structured abundance population model, especially with a
Beverton-Holt stock-recruitment relationship.
With this, we provide numerical estimates for two Chilean fisheries.

\subsection{An age class dynamical model}

We consider an age-structured abundance population model with a
possibly non linear stock-recruitment relationship, derived from
fish stock management (see~\cite{Quinn-Deriso:1999}, and also 
\cite{DeLara-Doyen-Guilbaud-Rochet_IJMS:2007} for more details). 
\medskip

Time is measured in years, and the time index $t \in \NN$
represents the beginning of year $t$ and of yearly period
$[t,t+1[$. Let $A \in \NN^*$ denote a maximum age, and $a \in \{
1, \ldots, A\}$ an age class index, all expressed in years. 
The state is the vector $\abundance=(\abundance\ind{a})_{a=1,\ldots, A}
\in \RR^A_+$, the \emph{abundances} at age: for $a=1,\ldots, A-1$,
$\abundance\ind{a}(t)$ is the number of individuals of age between
$a-1$ and $a$  at the beginning of yearly period $[t,t+1[$;
$\abundance\ind{A}(t)$ is the number of individuals of age greater
than $A-1$.
The control $\multiplier(t) $ is the  \emph{fishing effort multiplier},
 supposed to be  applied in the middle of period $[t,t+1[$.
The control dynamical model is 
\begin{equation*}
\abundance(t+1) = 
\dynamics \big(\abundance(t), \multiplier(t) \big), \quad
t= t_{0}, t_{0}+1, \ldots , \qquad \abundance(t_{0}) \mtext{
given,} \label{dynstock}
\end{equation*}
where the vector function $ \dynamics = \ds \left(
\dynamics\ind{a} \right)_{a=1, \hdots, A}$ is defined   for any $
\abundance \in \RR_+^A$ and $\multiplier \in \RR_+$ by
\begin{equation}
\left\{
\begin{array}{lcl}
\dynamics\ind{1}(\abundance,\multiplier) &=&
\varphi \big( \SSB(\abundance) \big) \, , \\[4mm]
\dynamics\ind{a}(\abundance,\multiplier) &=&
e^{-(\mortality\ind{a-1} + \multiplier \exploitation\ind{a-1})}
\abundance\ind{a-1},
\quad a = 2, \hdots, A-1 \, , \\[4mm]
\dynamics\ind{A}(\abundance,\multiplier) &=&
e^{-(\mortality\ind{A-1} + \multiplier \exploitation\ind{A-1})}
\abundance\ind{A-1} + \pi \times e^{-(\mortality\ind{A} +
\multiplier \exploitation\ind{A})} \abundance\ind{A} \, .
\end{array}
\right. \label{funcdyn}
\end{equation}
In the above formulas, 
$\mortality\ind{a}$ is the natural \emph{mortality rate} of
  individuals of age $a$, 
$\exploitation\ind{a}$ is the mortality rate of individuals
of age   $a$ due to   harvesting between $t$ and $t+1$, supposed to remain
  constant during   period  $[t,t+1[$
(the vector $(\exploitation\ind{a})_{a=1,\ldots,A}$ is termed the
\emph{exploitation pattern}), 
and the parameter $\pi \in \{0,1\}$ is related to the existence of a
so-called \emph{plus-group} (if we neglect the survivors older
than age $A$ then $\pi=0$, else $\pi=1$ and the last age class is
a plus group).
The function $\varphi$ describes a \emph{stock-recruitment
relationship}.
The \emph{spawning stock biomass} $\SSB$ is defined by
\begin{equation}
\SSB(\abundance) \defegal \sum_{a=1}^{A} \mature\ind{a}
\weight\ind{a} \abundance\ind{a} \, , \label{eq:SSB}
\end{equation}
that is summing the contributions of individuals to reproduction,
where $(\mature\ind{a})_{a=1, \ldots, A}$ are the
\emph{proportions of mature individuals} (some may be zero) at age
and $(\weight\ind{a})_{a=1, \ldots, A}$ are the \emph{weights at
age} (all positive). 

\subsection{An acceptable set reflecting conflicting preservation and
production objectives}

We shall consider an acceptable set $\obj$ which reflects conflicting
objectives of \emph{preservation} -- measured by the 
spawning stock biomass being high enough -- 
and of \emph{production}, measured by the following yield indicator. 

The exploitation is described by catch-at-age $\catch\ind{a}$ and
yield $\Yield$, defined for a given vector of
abundance $\abundance$ and a given control $\multiplier$ by the so
called \emph{Baranov catch equations}
\cite[p.~255-256]{Quinn-Deriso:1999}. The catches are the number
of individuals captured over the period $[t-1,t[$:
\begin{equation*}
\catch\ind{a} \bigl(\abundance,\multiplier\bigr)= \frac{
\multiplier \exploitation\ind{a}}{ \multiplier
\exploitation\ind{a} +\mortality\ind{a}} \left(
1-e^{-(\mortality\ind{a}+ \multiplier
\exploitation\ind{a})}\right) \abundance\ind{a} \, .
\label{eq:SSVPA_catch}
\end{equation*}
The production in term of biomass at the beginning of period
$[t,t+1[$ is then
\begin{equation}
\Yield\bigl(\abundance,\multiplier\bigr)=
 \sum_{a=1}^{A} \weight\ind{a}\, \catch\ind{a}(\abundance,\multiplier)
 \; .
\label{eq:productionVPA}
\end{equation}

We focus our analysis on the acceptable set 
\begin{equation}
\obj_{\yield}(y_{\mmin},\Blim) \defegal  \{(\abundance,\multiplier) 
\mid \Yield(\abundance,\multiplier ) \geq y_{\mmin},~
\SSB(\abundance) \ge \Blim   \} \, , 
 \end{equation}
where the yield function $\Yield$ is given by
\eqref{eq:productionVPA} and $\SSB$ by~\eqref{eq:SSB}. 
Contrarily to the ICES precautionary approach 
as analyzed in \cite{DeLara-Doyen-Guilbaud-Rochet_IJMS:2007},
we not only focus on \emph{preservation issues} 
($\SSB(\abundance) \ge \Blim $) but also on 
\emph{production issues} by asking for a minimal yield
($\Yield(\abundance,\multiplier ) \geq y_{\mmin}$).

\subsection{Monotonicity properties}

The  set $\obj_{\yield}(y_{\mmin},\biomass_{\mmin})$ 
is a production acceptable set. Indeed, on the one hand,
the yield $\Yield$ is increasing with respect both to the state and to the
control. 
On the other hand, the spawning stock biomass $\SSB$ is increasing with
respect to the state (and does not depend on the control).

The dynamics~\eqref{funcdyn} is a monotone bioeconomic one whenever 
the recruitment function $\varphi$ in \eqref{funcdyn} is non
decreasing. 

We now focus on the existence of equilibrium points.
As is classical, consider the following  proportions  of equilibrium
recruits which survive up to age $a$:
\begin{equation*} 
\left\{ \begin{array}{rcl}
s\ind{1}(\multiplier) &\defegal &1 \\
 s\ind{a}(\multiplier) &\defegal &\ds \exp \Bigl( - \bigl( \mortality\ind{1} +
\cdots + \mortality\ind{a-1} + \multiplier(\exploitation\ind{1} +
\cdots + \exploitation\ind{a-1}) \bigr) \Bigr) \, , \quad
a=2,\ldots,A-1 \\
 s\ind{A}(\multiplier) &\defegal &\ds \frac{1}{%
1 - \pi e^{-(\mortality\ind{A}+\multiplier \exploitation\ind{A})}
} \exp \Bigl( - \bigl( \mortality\ind{1} + \cdots +
\mortality\ind{A-1} + \multiplier(\exploitation\ind{1} + \cdots +
\exploitation\ind{A-1}) \bigr) \Bigr)
\end{array} \right.
\end{equation*}
Let also \(\spr(\multiplier) \defegal \sum_{a=1}^{A} \mature\ind{a}
\weight\ind{a} s\ind{a}(\multiplier)\) be the 
\emph{spawning per recruit} at equilibrium.
 When the recruitment function  is Beverton-Holt
$\varphi(\biomass)= \frac{\biomass}{\alpha+\beta \biomass}$ (which
includes the constant case, taking $\alpha=0$,
 and the linear case, taking $\beta=0$),
there exists an  equilibrium point for any control $\multiplier \geq 0$.
It is given by 
$\overline{\abundance}(\multiplier)=
(\overline{\abundance}\ind{a}(\multiplier))_{a=1,\ldots,A}$, where
\( \overline{\abundance}\ind{a}(\multiplier)=Z(\multiplier)
s\ind{a}(\multiplier) \) 
and 
\(  Z(\multiplier) = 
\max\left\{0,\frac{\spr(\multiplier)-\alpha}{\beta\spr(\multiplier)}\right\}
\) if $\beta >0$,
\(  Z(\multiplier) = 0 \) if $\beta =0$.

\if{
We know give conditions for the equilibrium point 
$\overline{\abundance}(\multiplier\llower)$
to be globally asymptotically stable on $ \SCS$.

\begin{lemma}
Assume that the stock-recruitment relationship $\varphi$
 is Beverton-Holt
$\varphi(\biomass)= \frac{\biomass}{\alpha+\beta \biomass}$ (allowing
the cases $\alpha=0$ or $\beta=0$). If 
\begin{equation}
\label{eq:Blim_one}
\left\{\begin{array}{ll}
\Blim \ge \overline{\Blim}({\multiplier\llower}) \defegal  
\frac{\sqrt{\frac{\alpha ~\max \limits_{a =1,\hdots,A}
      \mature\ind{a}\weight\ind{a}}{1- \max \limits_{a =1,\hdots,A}
      e^{-(\mortality\ind{a} + 
\multiplier\llower \exploitation\ind{a})} }} ~ -  
\alpha}{\beta} & \mbox{if }~\beta > 0\\[3mm]
\frac{\max \limits_{a =1,\hdots,A} \mature\ind{a}\weight\ind{a}}{\alpha}
+ \max \limits_{a =1,\hdots,A} e^{-(\mortality\ind{a} + 
\multiplier\llower \exploitation\ind{a})}  < 1 & \mbox{if }~\beta =0 \; ,
 \end{array} \right.
\end{equation}
the equilibrium point $\overline{\abundance}(\multiplier\llower)$
is globally asymptotically stable for the dynamics
$\dynamics(\cdot,\multiplier\llower)$ on $\SCS$. 
\label{lem:Blim_one}
\end{lemma}

\begin{proof}
The projection $\SCS$ of the above desirable state onto the
state space $\RR^A_+$ is included in the set 
\(
\obj(\Blim) \defegal \{\abundance ~\mid ~\SSB(\abundance) \ge \Blim\} 
\).
Since, for any $\abundance$, $\abundance'$ in $\obj(\Blim)$  we have
\[
\|\dynamics(\abundance,\multiplier\llower)-\dynamics(\abundance',\multiplier\llower)\|_1
\le \left( \max_{B \ge \Blim}|\varphi'(B)| 
\max_{a =1,\hdots,A} \mature\ind{a}\weight\ind{a} + 
\max_{a =1,\hdots,A} e^{-(\mortality\ind{a} +
\multiplier\llower \exploitation\ind{a})} \right) \|\abundance -
\abundance'\|_1 \; , 
\]
where $\|\abundance\|_1$ is the norm $\sum_{a=1}^A N\ind{a}$ in
$\RR^A_+$, 
a sufficient condition for the existence of a globally asymptotically
stable steady state for the dynamics
$\dynamics(\cdot,\multiplier\llower)$ on $\SCS$ is to have  
\[
 \max_{B \ge \Blim}|\varphi'(B)| \max_{a =1,\hdots,A} \mature\ind{a}\weight\ind{a} + \max_{a =1,\hdots,A} e^{-(\mortality\ind{a} +
\multiplier\llower \exploitation\ind{a})}  < 1 \; .
\]
When the function $\varphi$  is Beverton-Holt
$\varphi(\biomass)= \frac{\biomass}{\alpha+\beta \biomass}$,
$\varphi'$ is decreasing ($\varphi$ is concave), 
and the above condition is equivalent to~\eqref{eq:Blim_one}.
\end{proof}


}\fi

\subsection{Minimal viable production  issues}

The following statement establishes 
\emph{maximum sustainable thresholds} for
the indicators $\SSB$ and $\Yield$. 
It is an application of Corollary~\ref{corollary1}.

 \begin{proposition}
\label{prop:mvp}
Assume that the stock-recruitment relationship $\varphi$
 is Beverton-Holt
$\varphi(\biomass)= \frac{\biomass}{\alpha+\beta \biomass}$ (allowing the cases $\alpha=0$ or $\beta=0$),  that the fishing effort
$\multiplier$ is bounded from below and above by
$0 \leq \multiplier\llower \leq \multiplier \leq \multiplier\upper$. If 
\begin{equation}
\label{lipschitz}
\phi_\dynamics(\multiplier\llower)\defegal \varphi'\left(\SSB(\overline{\abundance}(\multiplier\llower))\right) \max_{a =1,\hdots,A} \mature\ind{a}\weight\ind{a} + \max_{a =1,\hdots,A} e^{-(\mortality\ind{a} +
\multiplier\llower \exploitation\ind{a})}  < 1
\end{equation}
then, ensuring a minimal
viable production and spawning stock biomass requires that the 
production and preservation thresholds $y_{\mmin}$ and $\Blim$
be not too high:
\begin{equation*}\label{prop:NC}
 \left. \begin{array}{r}
y_{\mmin} >
\Yield(\overline\abundance(\multiplier\llower),\multiplier\upper)
\\[3mm]
\mbox{or }~ \Blim > \SSB(\overline{\abundance}(\multiplier\llower))
\end{array}\right\} \Rightarrow 
\VV\big( \dynamics,\obj_{\yield}(y_{\mmin},\Blim) \big)
= \emptyset \; .
\end{equation*}
\label{prop:main_fish}
\end{proposition}

\begin{proof}
In order to apply Corollary ~\ref{corollary1}, 
let us prove that, for $\Blim >
\SSB(\overline{\abundance}(\multiplier\llower))$, one has the following
property  
\begin{equation}\label{def_cte}
\|\dynamics(\abundance,\multiplier\llower)-\overline{\abundance}(\multiplier\llower)\|_1
\le \phi_\dynamics(\multiplier\llower) \|\abundance -
\overline{\abundance}(\multiplier\llower)\|_1 \; , 
\end{equation}
for all $\abundance$ in   $\SCS$ (projection on $\RR^A_+$  of the
acceptable set $\obj_{\yield}(y_{\mmin},\Blim)$)  where
$\|\abundance\|_1$ is the norm $\sum_{a=1}^A |N\ind{a}|$ in $\RR^A$.
For any $\abundance$ one has
\begin{equation}\label{difN}
\|\dynamics(\abundance,\multiplier\llower)-\overline{\abundance}(\multiplier\llower)\|_1\le \left|\varphi(\SSB(\abundance))-\varphi(\SSB(\overline{\abundance}(\multiplier\llower))\right| + \sum_{a=1}^A e^{-(\mortality\ind{a} +
\multiplier\llower \exploitation\ind{a})}
|\abundance_a-\overline{\abundance}(\multiplier\llower)_a| \; .
\end{equation}
If $\abundance \in \SCS$ then $\SSB(\abundance) \ge \Blim >
\SSB(\overline{\abundance}(\multiplier\llower))$ and therefore
one obtains
\[
\left|\varphi(\SSB(\abundance))-\varphi(\SSB(\overline{\abundance}(\multiplier\llower))\right|
\le \max_{B \in
  \left[\SSB(\overline{\abundance}(\multiplier\llower)),\SSB(\abundance)\right]} |\varphi'(B)| \sum_{a=1}^A \mature\ind{a}\weight\ind{a}|\abundance_a-\overline{\abundance}(\multiplier\llower)_a| \; .
\] 
By the concavity of $\varphi$ (which implies that $\varphi'$ is
decreasing), this gives
\begin{eqnarray*}
\left|\varphi(\SSB(\abundance))-\varphi(\SSB(\overline{\abundance}(\multiplier\llower))\right|
& \leq & 
\left|\varphi'(\SSB(\overline{\abundance}(\multiplier\llower)))\right| \sum_{a=1}^A
\mature\ind{a}\weight\ind{a}|\abundance_a-\overline{\abundance}(\multiplier\llower)_a|
\\ & \leq & 
\left(\varphi'\left(\SSB(\overline{\abundance}(\multiplier\llower))\right) \max_{a =1,\hdots,A} \mature\ind{a}\weight\ind{a}\right) \|\abundance -
\overline{\abundance}(\multiplier\llower)\|_1 \; .
\end{eqnarray*}
The above inequality together with \eqref{difN} and the definition of
\eqref{lipschitz} make it possible to obtain \eqref{def_cte} and, then, the
condition 
 \eqref{condition_steady}  of Corollary ~\ref{corollary1}.
\end{proof}

The above result can be interpreted as follows.
\begin{itemize}
\item There is no vector of abundance which allows to obtain,
  starting  from it, catches greater than the
\emph{maximal production threshold}
  $\Yield(\overline\abundance(\multiplier\llower),\multiplier\upper)$,
  during all the periods.

\item  Starting from any vector of abundance, 
whatever the harvest, the minimum level of spawning stock biomass
($\SSB$) observed during all the  periods will be lower than 
(or equal to) the 
\emph{maximal preservation threshold}
$\SSB(\overline{\abundance}(\multiplier\llower))$. 
\end{itemize}

\subsection{Numerical applications to Chilean fisheries}
\label{sec:numerical}

We provide numerical estimates obtained for the species Chilean sea bass
(\emph{Dissostichus eleginoides}), harvested in the south of Chile, and
Alfonsino (\emph{Beryx splendens}), harvested in the Juan Fernández
archipelago. 
The dynamics of the Chilean sea bass can be described by the
model~\eqref{funcdyn} with a Beverton-Holt 
stock-recruitment relationship $\varphi$. 
For the Alfonsino, females and males are distinguished, each following 
a dynamics~\eqref{funcdyn} with a Beverton-Holt 
stock-recruitment relationship $\varphi$. 
Thus, for this species, the state is the abundances at age for
females and males and the resulting dynamics  is a monotone bioeconomic
one. For both species,  the mortality is supposed to be  the same at all
ages, and will be denoted by $\mortality$.
In Tables~\ref{table:data}, \ref{table:age_bacalao},
and~\ref{table:age_alfonsino} in the Appendix,
we present  numerical data for 
both species, provided by the \emph{Centro de Estudios Pesqueros -
  Chile (CEPES)}. 
 
Table~\ref{table:mt} sums up the maximal production and preservation
thresholds obtained from Proposition~\ref{prop:main_fish} for both
species and the values of $\phi_\dynamics(\multiplier\llower)$
 defined by  \eqref{lipschitz}.

 \begin{table}[h]
\begin{center}
{\scriptsize 
\begin{tabular}{|lcrr|} 
\hline
 Definition & Notation& \mbox{\it Chilean sea bass}  & \mbox{\it Alfonsino} \\\hline\hline
 Maximal threshold for a sustainable catch (tons) & $ \Yield(\overline\abundance(\multiplier\llower),\multiplier\upper) $ &15~166 & 16~158\\\hline
Maximal threshold for a sustainable $\SSB$ (tons) &$\SSB(\overline{\abundance}(\multiplier\llower))$ &56~521 & 52~373\\\hline
 Constant defined by \eqref{lipschitz}  &$\phi_\dynamics(\multiplier\llower)$ & 0.852 & 0.818 \\\hline
\end{tabular}
\caption{Maximal sustainable thresholds for {\it Chilean sea bass}
and for {\it Alfonsino}.
\label{table:mt}}}
\end{center}
\end{table}


\subsection{Chilean sea bass}

\begin{center}
\begin{figure}[h]
\begin{center}
\includegraphics[height=7cm,width=13cm]{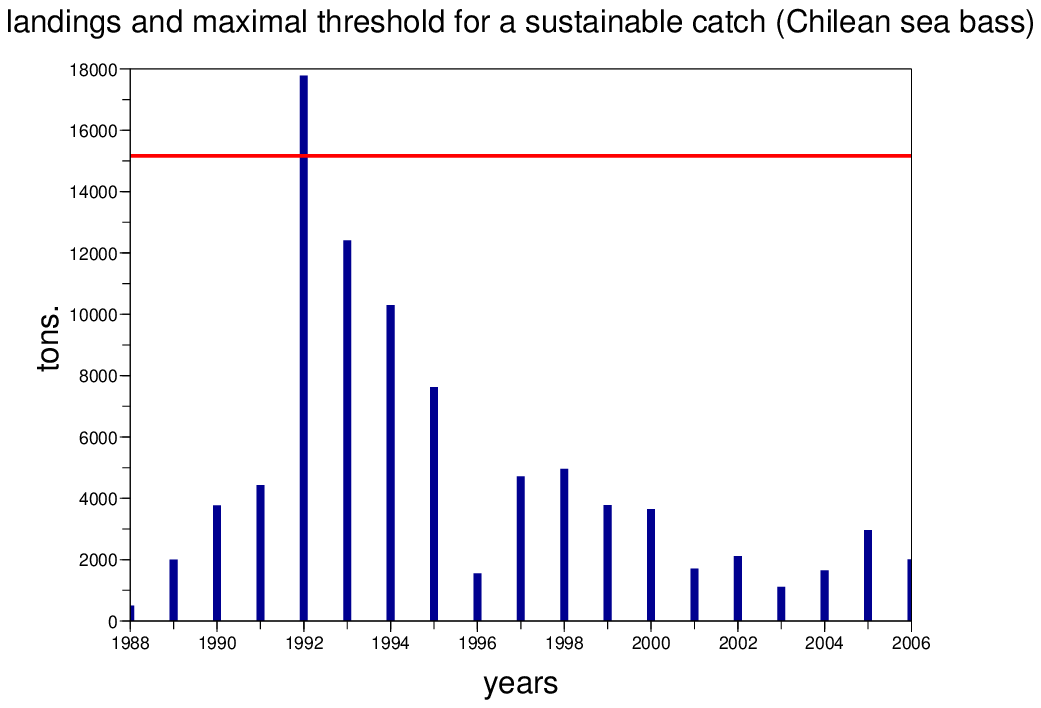}
\caption{Chilean sea bass: landings (1988-2006) [tons] and $\Yield(\overline\abundance(\multiplier\llower),\multiplier\upper)$.}
\label{Fig:bacalao_ymin}
\end{center}
\end{figure}
\end{center}

Figure~\ref{Fig:bacalao_ymin} displays the Chilean sea bass  landings,
between 1988 and 2006.
The horizontal line represents the maximal threshold
$\Yield(\overline\abundance(\multiplier\llower),\multiplier\upper)$.
Hence, it may be seen that the catches obtained in 1992 were not
sustainable: even if the species were abundant, such landings could not
be maintained forever.

\begin{center}
\begin{figure}[h]
\begin{center}
\includegraphics[height=7cm,width=13cm]{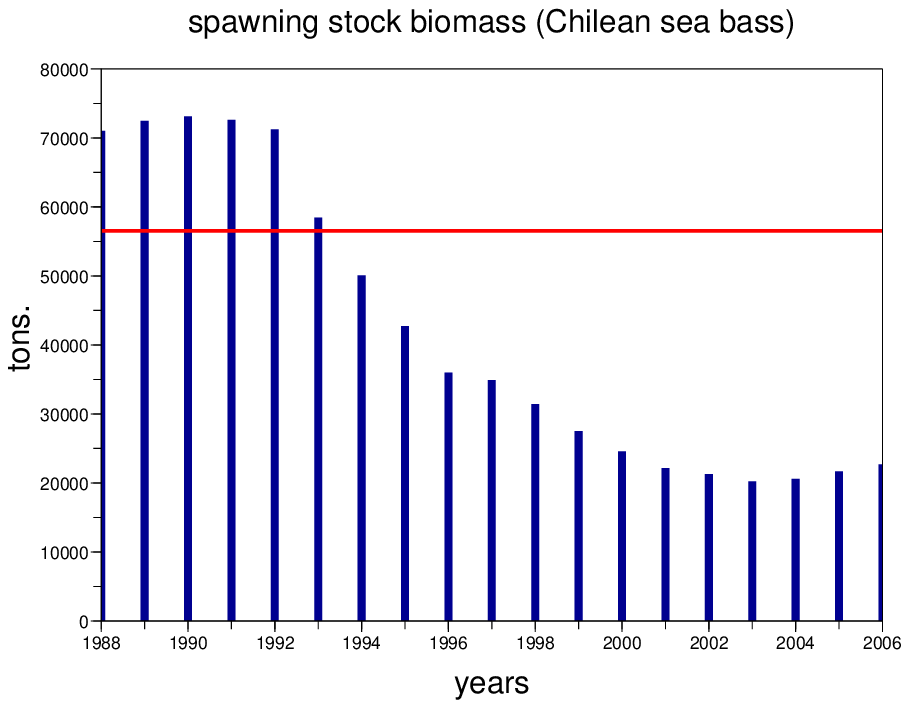}
\caption{Chilean sea bass: SSB (1988-2006) [tons] and $\SSB(\overline \abundance(\multiplier\llower))$.}
\label{Fig:bacalao_SSB}
\end{center}
\end{figure}
\end{center}

Figure~\ref{Fig:bacalao_SSB} displays the Chilean sea bass  
spawning stock biomass ($\SSB$), between 1988 and 2006.
The horizontal line represents the maximal threshold 
$\SSB(\overline \abundance (\multiplier\llower))$.  
The $\SSB$ observed during the first six years could not have been
sustained forever.

\subsection{Alfonsino}

\begin{center}
\begin{figure}[hhhh]
\begin{center}
\includegraphics[height=7cm,width=13cm]{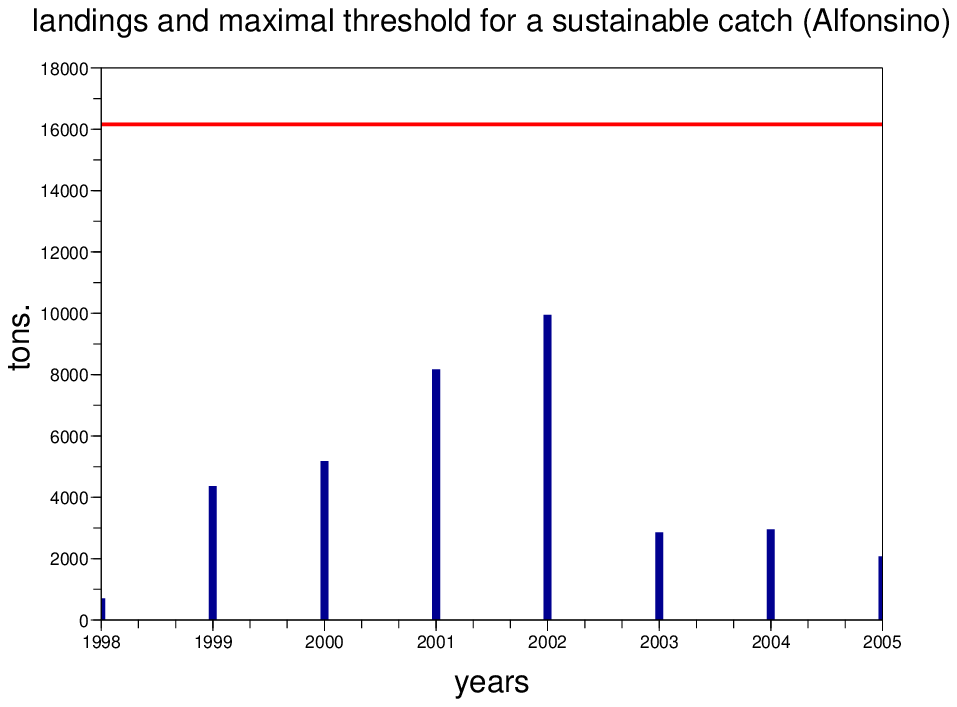}
\caption{Alfonsino: landings (1998-2005) [tons] and $\Yield(\overline\abundance(\multiplier\llower),\multiplier\upper)$.}
\label{Fig:alfonsino_ymin}
\end{center}
\end{figure}
\end{center}

\begin{center}
\begin{figure}[hhhh]
\begin{center}
\includegraphics[height=7cm,width=13cm]{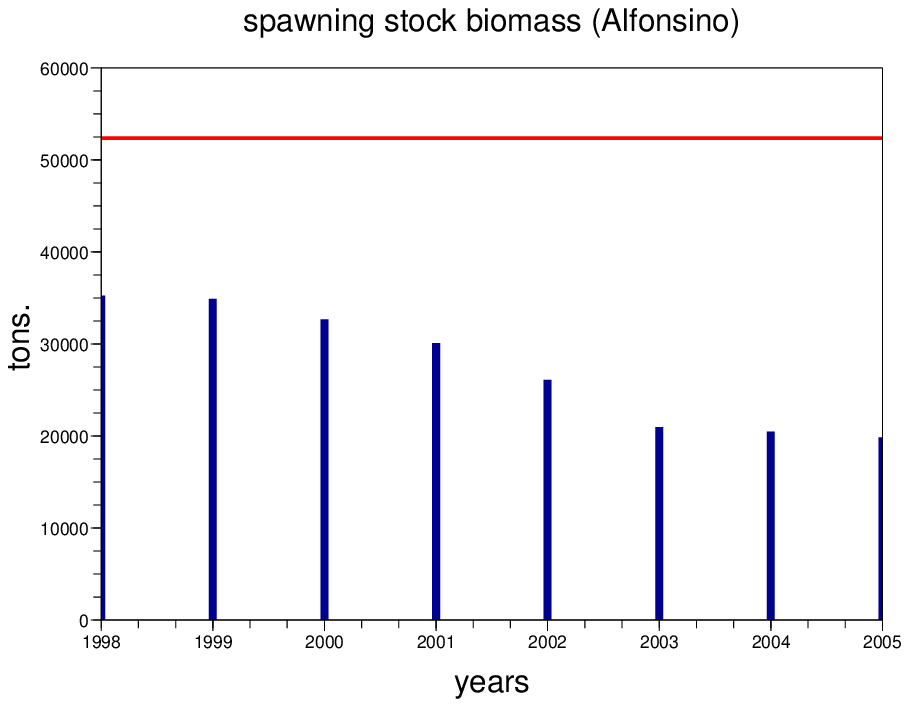}
\caption{Alfonsino: SSB (1998-2005) [tons] and $\SSB(\overline \abundance(\multiplier\llower))$.}
\label{Fig:alfonsino_SSB}
\end{center}
\end{figure}
\end{center}

For the Alfonsino, Figures~\ref{Fig:alfonsino_ymin}
and~\ref{Fig:alfonsino_SSB}, both spawning stock biomasses and 
landings are below the maximal threshold. 
Thus, we cannot conclude that these levels indicate a non viable fishery
management. 

\section{Conclusion}

Some monospecies age class models, as well as specific multi-species
models (with so-called technical interactions), exhibit useful 
monotonicity properties.
We have shown how these latter may help providing estimates of the 
viability kernel for so-called 
production and preservation acceptable sets.
When the acceptable set is defined by inequalities requirements given by
indicator functions, we provide conditions on the corresponding
thresholds to test whether the viability kernel is empty or not.

This theoretical framework is applied to fishery management analysis. 
We obtain upper bounds for production which are interesting for managers
in that they only depend on the model parameters, 
and not on the current stocks.
Our formulas for so-called maximal sustainable thresholds give sensible
values: Chilean sea bass data violate these bounds, 
while Alfonsino data are within.

We have thus provided a general method to analyze up to what points can 
conflicting production and preservation objectives be sustainably
achieved for a class of models including monospecies age class 
and multi-species with technical interactions.
\bigskip

\begin{acknowledgements}

This paper was prepared within the MIFIMA (Mathematics, 
  Informatics and Fisheries Management) 
international research network.
We thank CNRS, INRIA and the French Ministry of Foreign Affairs 
for their funding and support through the regional cooperation program
STIC--AmSud.
We also thank  CONICYT (Chile) for its support  through projects
 ECOS-CONICYT number C07E03, STIC--AmSud, FONDECYT N 1070297 (H. Ramírez C.), FONDECYT N 1080173 (P. Gajardo), and Fondo Basal, Centro de Modelamiento Matematico, U. de Chile. 
 Finally, authors are indebted to Alejandro Zuleta and Pedro Rubilar, from CEPES (Centro de Estudios Pesqueros), Chile, for their contribution to Section 4 of the current paper.

\end{acknowledgements}

\bibliographystyle{spbasic}      



\newcommand{\noopsort}[1]{} \ifx\undefined\allcaps\def\allcaps#1{#1}\fi

\newpage 
\appendix 

\section{Appendix: model data}

 \begin{table}[h]
\begin{center}
{\scriptsize 
\begin{tabular}{|lccc|} 
\hline
 Definition & Notation& \mbox{\it Chilean sea bass}  & \mbox{\it Alfonsino} \\\hline\hline
 Maximum age &$A$ & 36 & 40 (20 male \& 20 female) \\\hline
 Beverton-Holt parameter (adimensional)  &$\alpha$& $1.4 \cdot 10^{-3}$&$9.16 \cdot 10^{-5}$  \\\hline
  Beverton-Holt parameter (gr) &$\beta$& $4.65 \cdot 10^{-7}$ & $7.46 \cdot 10^{-8}$ \\\hline
 Natural mortality &$M$& 0.16 & 0.2 \\\hline
 Presence of plus group &$\pi$  & $1$ & $1$  \\\hline
Lower limit for fishing effort multiplier &$\multiplier\llower$  & $0$ & $0$ \\\hline
Upper limit for fishing effort multiplier &$\multiplier\upper$  & $0.39$ & $0.885$ \\\hline
\end{tabular}
\caption{Parameter definitions and values for two case studies.
}
\label{table:data} 
}
\end{center}
\end{table}

 \begin{table}[h]
\begin{center}
{\scriptsize 
\begin{tabular}{|c|c|c|c|} 
\hline
Age $a$ &Mean weight-at-age (gr) $~(\weight\ind{a})_a$& Maturity ogive $~(\mature\ind{a})_a$ &  Fishing mortality-at-age $(\exploitation_a)_a$  \\\hline\hline
1	&3&	0&	0.0005\\\hline2	&77&	0&	0.0013\\\hline3	&326&	0&	0.0051\\\hline4	&809&	0&	0.0183\\\hline5	&1547&	0&	0.0494\\\hline6	&2536&	0&	0.1080\\\hline7	&3753&	0&	0.2067\\\hline8	&5169&	0.1& 	0.3467\\\hline9	&6748&	0.2&	0.5277\\\hline10	&8454&	0.3&	0.7127\\\hline11	&10253&	0.5&	0.8675\\\hline12	&12113&	0.7&	0.9611\\\hline13	&14006&	0.8&	1.0000\\\hline14	&15905&	0.9&	0.9831\\\hline15	&17791&	1&	0.9302\\\hline16	&19646&	1&	0.8661\\\hline17&	21456&	1&	0.7933\\\hline18	&23209&	1&	0.7254\\\hline19	&24897&	1&	0.6614\\\hline20	&26514&	1&	0.6040\\\hline21	&28056&	1&	0.5537\\\hline22	&29520&	1&	0.5091\\\hline23	&30906&	1&	0.4701\\\hline24	&32213&	1&	0.4434\\\hline25	&33441&	1&	0.4115\\\hline26	&34594&	1&	0.3867\\\hline27	&35673&	1&	0.3651\\\hline28	&36681&	1&	0.3464\\\hline29	&37621&	1&	0.3341\\\hline30	&38496&	1&	0.3206\\\hline31	&39309&	1&	0.3089\\\hline32	&40064&	1&	0.2986\\\hline33	&40763&	1&	0.2896\\\hline34	&41411&	1&	0.2817\\\hline35	&42011&	1&	0.2747\\\hline36	&45409&	1&	0.2408\\\hline
\end{tabular}
\caption{Parameters at age for the Chilean sea bass.
}
\label{table:age_bacalao}
}
\end{center}
\end{table}
 
  \begin{table}[h]
\begin{center}
{\scriptsize 
\begin{tabular}{|c|c|c|c|}
\hline
Age $a$ (male \& female) &Mean weight-at-age (gr) $~(\weight\ind{a})_a$ & Maturity ogive $~(\mature\ind{a})_a$ &  Fishing mortality-at-age $(\exploitation_a)_a$  \\\hline\hline
1&	138 $\slash$140 &	0	&0,006	$\slash$0,006\\\hline2&	252$\slash$	256	&0&	0,026$\slash$	0,024\\\hline3&	397$\slash$	406	&0&	0,083$\slash$	0,075\\\hline4&	565$\slash$	583	&0&	0,217$\slash$	0,192\\\hline5&	752$\slash$	783	&0&	0,441	$\slash$0,390\\\hline6&	951$\slash$	999	&1&	0,691$\slash$	0,619\\\hline7&	1156$\slash$	1227	 &1&	0,883$\slash$	0,809\\\hline8&	1364$\slash$	1462 &	1&	0,997	$\slash$0,932\\\hline9&	1571$\slash$	1699	 &1&	1,000$\slash$	1,000\\\hline10&	1774$\slash$	1935	 &1&	0,564$\slash$	0,884\\\hline11&	1970$\slash$	2168	 &1&	0,229$\slash$	0,486\\\hline12&	2158$\slash$	2395	 &1&	0,098$\slash$	0,233\\\hline13&	2337$\slash$	2614	 &1&	0,045$\slash$	0,118\\\hline14&	2506$\slash$	2825	 &1&	0,022$\slash$	0,064\\\hline15&	2664$\slash$	3026	 &1&	0,012$\slash$	0,037\\\hline16&	2812$\slash$	3216	 &1&	0,006$\slash$	0,022\\\hline17&	2950$\slash$	3396	 &1&	0,004$\slash$	0,014\\\hline18&	3078$\slash$	3565	 &1&	0,002$\slash$	0,010\\\hline19&	3196$\slash$	3724	 &1&	0,001$\slash$	0,007\\\hline20&	3304$\slash$	3872	 &1&	0,001$\slash$	0,005\\\hline
\end{tabular}
\caption{Parameters at age for the Alfonsino (male \& female). Maturity ogive parameters are the same for males and females.
}
\label{table:age_alfonsino}
}
\end{center}
\end{table}

\end{document}